\pdfoutput=1
\documentclass[12pt]{amsart}
\usepackage[leqno]{amsmath}
\usepackage{hyperref}
\usepackage{enumerate}
\usepackage{mathrsfs}
\usepackage{pdfpages}

\makeatletter
\newcommand{\printaddress}{\@setaddresses}
\makeatother

\newcommand{\tuple}[1]{\ensuremath{\left \langle #1 \right \rangle }}

\theoremstyle{plain}
  \newtheorem*{claim}{CLAIM}

\numberwithin{equation}{section}

\title[On the impossibility of using analogue machines]{On the impossibility of using analogue machines to calculate non-computable functions}
\author{R.O. Gandy\textsuperscript{\dag} \\ September 1993}
\thanks{\textsuperscript{\dag}Typeset by Aran Nayebi on August 27, 2013. A.N. is grateful to S. Barry Cooper and Philip Welch for providing a photocopy of Gandy's original handwritten manuscript (attached at the end of this document), as well as Solomon Feferman for suggesting to typeset it and his support.}
\address{9 Squitchey Lane, Oxford OX2 7LD}
\begin{document}

\maketitle
\section*{Introduction}
A number of examples have been given of physical systems (both classical and quantum mechanical) which when provided with a (continuously variable) computable input will give a non-computable output. It has been suggested that these systems might allow one to design analogue machines which would calculate the values of some number-theoretic non-computable function. Analysis of the examples show that the suggestion is wrong. In \S 4 I claim that given a reasonable definition of analogue machine it will always be wrong. The claim is to be read not so much as a dogmatic assertion, but rather as a challenge. \newline
\indent In \S's 1 and 2 I discuss analogue machines, and lay down some conditions which I believe they must satisfy. In \S 3 I discuss the particular forms which a paradigm undecidable problem (or non-computable function) may take. In \S's 5 and 6 I justify any claim for two particular examples lying within the range of classical physics, and in \S 7 I justify it for two (closely connected) examples from quantum mechanics, and discuss, very briefly, other possible quantum mechanical situations. \S 8 contains various remarks and comments. In \S 9 I consider the suggestion made by Penrose that a (future) theory of quantum gravity may predict non-locally-determined, and perhaps non-computable patterns of growth for microsopic structures. My conclusion is that such a theory will have to have non-computability built into it.

\section{Analogue machines}
By a continuously variable quantity (`CVQ') I mean a physical quantity which is represented mathematically by a point in a metric space - e.g., by a real number, or a point of Hilbert space. This is not put forward as an exact definition, but as an indication of how I use the term. For CVQ's very natural definitions of `computable' have been given in Pour-El \& Richards (1989); I shall to this book as CAP. Roughly speaking `$x$ is computable' means that $x$ is the limit of a sequence of finitely presented approximations and a modulus of convergence for the sequence can be computed. \newline
\indent In the theoretical treatment of a physical device the CVQ's have exact values, and no bound is place, \underline{a priori} on their magnitude. But when such a device is to be used as an analogue machine to perform some calculation then there will be an upper limit $x$ on the size of a CVQ (an electric circuit will melt if the current is too large) and a lower limit $\epsilon$ on the accuracy with which it can be controlled or measured. Numerical values for $x$ and $\epsilon$ depend, of course, on the choice of units for the particular CVQ considered, but the ratio $x/\epsilon$ does not; so we define the \underline{precision ratio} (PR for short) of the CVQ to be $x/\epsilon$. In the theory of the machines there may be different variables having the same physical dimension; these are to be counted as distinct CVQ's and may have different precision rations. The `independent' variable time is also a CVQ and has a PR; when an analogue machine is to be \underline{used} some limit must be placed on its run-time. \newline
\indent We are concerned with matters of principle rather than of practice, so although a given analogue machine will have definition precision ratios, we do not place any bound on the PR's that may be attained by some machine. \newline
\indent We are primarily - sometimes only - concerned with those CVQ's which are inputs and outputs of the machine. We shall be interested in cases where these may be continuously controlled or continuously recorded functions; in such cases the relevant $x$ and $\epsilon$ will be given by some norm for the functions. Most usually the uniform norm will be appropriate, but some machines one might want to say, the $L^2$ norm.\newline
\indent Even discretely varying quantities such as natural numbers have precision ratios attached to them; perfect accuracy (say $\epsilon < 1/2$) may be attainable, but still there is a bound on the size: one cannot place more than $N$ balls in a given box nor record more than $N$ events with a given geiger counter. In particular if an analogue machine incorporates a battery of digital computers then a PR (which depends both on the programme used and on the hardware) can be assigned to each of them; note that it does not depend on the placing of the decimal point. \newline
\indent In what follows we shall be concerned with the orders of magnitude of PR's rather than with precise values or upper bounds.

\section{Specification of analogue machines}
A specification for an analogue machine is a finite list of instructions which would, in principle, enable a technician or engineer to construct it; descriptions of the apparatus used in a (published) account of an experiment, do, although greatly abbreviated, have this form. If the correct operation of the machines requires particular precision ratios for certain quantities, then the instructions will specify tolerances for certain components\footnote{When A.M. Turing was building his speech encoder (`Delilah') he found that if it was to work, some of the components had to have a tighter than usual tolerance on their values; these were more expensive than the standard components and - at least in the case of resistances - had a gold spot to indicate that they were accurate to within (I think) 1\%.}. For example a machine might require a cam whose ideal shape ideal shape would be given by $r = f(\theta)$ where $f$ is some mathematical function. Then the instructions would indicate how the function $f$ could be computed (e.g., $f(\theta) = 2 + \sin^2{\theta}$ cms for $0 \le \theta \le 360^{\circ}$) and give a permitted tolerance (e.g., $\pm 10^{-3}$ cms). Tolerances can be given as precision ratios (${3.10}^3$ in the example). A specification will determine either explicitly or implicitly the PR's in the quantities (including outputs and inputs) occurring in the machine.

\section{Undecidable problems}
In the examples known to me it is proposed that there might be an analogue machine which with input $j(\in \mathbb{N})$ would output `Yes' or `No' to questions of the form ? $j \in A$? where $A$ is some standard recursively enumerable non-recursive set - for example the set which represents the halting problem. I shall only consider proposed machines of this kind. I describe two ways of representing the set $A$.
\subsection{}
There is a total computable function $a: \mathbb{N}\to\mathbb{N}$ which enumerates $A$ without repetitions. (This is the notation used throughout CAP). \newline
\indent The \underline{waiting-time} function $\nu$ is defined by
\begin{equation}
\nu(j) \simeq \mu n. \mbox{     }a(n) = j.
\end{equation}
This is a partial recursive function whose domain is $A$ and which is not bounded by any total computable function. For any particular analogue machine there is an upper bound $J$ on the inputs it can accept. I define
\begin{equation}
\beta(J) = \operatorname{Max}\{\nu(j): j < J \mbox{           } \& \mbox{           } j \in A\}
\end{equation}
(with $\operatorname{Max} \emptyset = 0$). This is a total function which is not computable; indeed it eventually majorises every computable function.
\subsection{}
There is a polynomial $P_A(y, \vec{x})$ such that
\begin{equation} \label{poly}
j \in A \leftrightarrow (\exists\vec{m}) \mbox{     }P_A(j, \vec{m}) = 0,
\end{equation}
where the variables of $\vec{m}$ $(= m_1, m_2, \ldots, m_k)$ range over the natural numbers. \newline
\indent In this case we define
\begin{equation}
\nu(j) \simeq (\mu n)(\exists \vec{m} < n)\mbox{     } P_A(j, \vec{m}) = 0,
\end{equation}
and
\begin{equation}
\beta(J) = \operatorname{Max}\{\nu(j): j \in A \mbox{           } \& \mbox{           }j < J\}.
\end{equation}
Then $\nu$ and $\beta$ have the same properties as in 3.1. Observe that, if $j \in A$, then
\begin{equation}\label{3.2.6}
\forall\vec{m} < \nu(j) \mbox{    }P_A(j, \vec{m}) \ne 0.
\end{equation}
\indent Various explicit definitions of suitable polynomials have been given. For each of these, if $P_A(j, \vec{m}) = 0$ then at least one of the $m_i$ encodes a \underline{particular} sequence which lists the first so many values of some recursive function. So, taking $i = 1$, we may suppose that
\begin{equation}\label{3.2.7}
P_A(j, \vec{m}) = 0 \textrm{ and }P_A(j, m', m_2, \ldots, m_k) \ne 0
\end{equation}
where $|m' - m| = 1$.

\section{The Claim}
Since a given machine cannot handle numbers greater than some bound we consider a given $J$ and the questions ?$j \in A$? for $j < J$. Now I make the following
\begin{claim}
Let $J$ be given. Then one cannot design an analogue machine (whose behaviour is governed by standard physical laws) which will give correct answers to all the questions ?$j\in A$? for $j < J$ unless one knows a bound $\beta$ for $\beta(J)$.
\end{claim}
\indent I call this a claim rather than a conjecture because I do not think one could prove it unless one placed severe restrictions on the notion of `analogue machine', and this I do not wish to do.\footnote{Pour-El in her (1974) gives a definition (based on differential analysers) of `General Purpose Analogue Computers' and characterizes the class of continuous functions which they can generate. She is not concerned with questions of precision, but I believe that the methods used in \S 5 and \S 6 can be applied to justify my claim for all machines of the type she considers.} But I believe that if someone proposes an analogue machine for settling ?$j \in A$? for $j < J$ then it can be shown that either they have (surreptiously?) made use of a bound for $\beta(j)$, or that not all the given answers will be correct. To illustrate the significance of the wording of the claim, suppose (what is quite plausible) that someone proves that $j \not\in A$ for all $j < J = 10$; then he can design a machine which always outputs `NO' for $j < J$. But, because of his proof he does in fact know that $\beta(J) = 0$. \newline
\indent Of course if one knows a $B$ as above then one does not need an analogue machines to settle ?$j \in A$? One simple computes $a(n)$ (as in 3.1) on $P_A(j,m_1,\ldots,m_k)$ (as in 3.2) for all $n < B$ or for all $m_1,\ldots,m_k < B$.

\section{First example (see CAP pp 51-53)}
\indent Let
\begin{equation} \label{1}
\phi(x) = \begin{cases}e^{-\frac{x^2}{1-x^2}}&\mbox{for } |x| \le 1 \\
0 & \mbox{for } |x| \ge 1. \end{cases} 
\end{equation}
$\phi$ is an infinitely differentiable function ($\in C^{\infty}$) though it is not analytic. Let
\begin{equation} \label{2}
\psi_n(x) = 4^{-a(n)}\phi\left(2^{-(n+a(n)+2)}\left(x-2^{-a(n)}\right)\right),
\end{equation}
where $a$ is as in \S 3.1. The graph of $\psi_n(x)$ is a blip of height $4^{-a(n)}$ centred on $2^{-a(n)}$, and having a width of $2^{-(n+a(n)+1)}$. If $m\ne n$ then the supports of $\psi_m$, $\psi_n$ do not intersect. Set
\begin{equation} \label{3}
f'(x) = \sum_{n = 0}^{\infty}\psi_n(x);
\end{equation}
$f'$ has a continuous but unbounded derivative, and $f'(x) = 0$ for $x > 5/4$. Since
\begin{equation} \label{4}
f'(2^{-j}) = \begin{cases}4^{-j}&\mbox{if } j \in A, \\
0 & \mbox{if } j \not\in A, \end{cases} 
\end{equation}
$f'$ is not a computable function. \newline
\indent Let
\begin{equation} \label{5}
\Phi_n(x) = \int_{0}^{x}\psi_n(x)dx.
\end{equation}
The graph of $\Phi_n$ is a smoothed out step function with initial value 0 (at $x = 0$) and a final value lying between 0 and $2^{-n}$.  \newline
\indent Now take
\begin{equation*}
f(x) = \sum_{n=0}^{\infty}\Phi_n(x).
\end{equation*}
$f$ is a computable function and its derivative is indeed the $f'$ given by \eqref{3}. Note that $||f||$, the uniform norm of $f$, is less than 2. To settle $j \in A$ the idea is to feed $f$ into an (analogue) differentiator, and then to observe whether the output $f'(x)$ is zero or not at $x = 2^{-j}$. For definiteness let us suppose that we control the current $i_1$ in a circuit $C_1$ inductively to a passive circuit $C_2$ and observe whether the current $i_2$ in $C_2$ is zero or not at time $2^{-j}$. The claim for this machine is justified on two counts.
\subsection{}
Because of the narrowness of the blip $\psi_n$, the measurement of the time $2^{-j}$ at which $i_2$ is observed must have, for $j \in A$, a precision ratio of order $2^{-\nu(j)}$ if the observed value of $i_2$ is to be different from zero.
\subsection{}
For $j \in A$, let
\begin{equation*}
f_j(x) = f(x) - \Phi_{\nu(j)}(x).
\end{equation*}
Then $f'_{j}(2^{-j}) = 0$. So if the machine is to give the answer YES for this $j$, then $i_1$ must satisfy
\begin{equation*}
|i_1(t) - f(t)| < \Phi_{\nu(j)}(t) \le 2^{-\nu(j)}.
\end{equation*}
So unless the precision ratio for the uniform norm of $i_1$ is better than $2^{\beta(J)}$ the machine will give wrong answers for some $j < J$.
\subsection{}
Thus to design a machine which will give correct answers for all $j < J$ we need to know $\beta(J)$.

\section{Second Example}
In their (1991) Doria \& Costa showed how a function defined in Richardson (1968) could theoretically be used in the construction (based solely on classical dynamics) of a device which would settle questions of the form ?$j \in A$?. They write\newline
\indent `Our example is intended to be seen as a \underline{Gedanken} experiment, as we do not wish to consider at the moment the certainly formidable question of its implementation.' \newline
\indent I shall show that its implementation by an analogue machine requires knowledge of a bound for $\beta(J)$.
\subsection{}
Let $k \ge 1$ be given and let $\mathscr{L}$ be the class of all real-valued functions of $k+1$ or fewer real variables which can be get by composition from the following initial functions:
\begin{enumerate}[(i)]
\item + and $\times$;
\item $\sin$;
\item projection functions $\lambda\vec{x}.x_i$;
\item constant functions $\lambda\vec{x}.c$, where $c$ is either $\pi$ or a rational number.
\end{enumerate}
\indent\indent Let $P_A$ be the polynomial of \S 2 \eqref{poly}. Richardson shows how one can define a function $F(u,x_1,\ldots,x_k)$ in $\mathscr{L}$ having the following properties.
\begin{enumerate}
\item $F$ is an even function of each of the $x_i$.
\item $F(u, x_1,\ldots,x_k) \ge 0$
\item $F(j,x_1,\ldots,x_k) > 1$ if $j \not\in A$.
\item If $F(j, x_1,\ldots,x_k) \le 1$ then $P_A(j,\tuple{x_1^2},\ldots,\tuple{x_k^2}) = 0$ and $F(j,\tuple{x_1^2},\ldots,\tuple{x_k^2}) = 0$ where $\tuple{x_i^2}$ denotes the natural number nearest to $x_i^2$. Hence in this case $j \in A$.
\item To calculate $F(j, x_1,\ldots,x_k)$ it is necessary first to calculate $P_A(j, x_1^2, \ldots, x_k^2)$
\end{enumerate}
\subsection{}
Let $\rho$ either be the function $\phi$ of \S 5, or be given by
\begin{equation*}
\rho(x) = \frac{1}{2}(|x-1|-(x-1)).
\end{equation*}
In either case $\rho(x) = 0$ for $x \ge 1$ and $\rho(0) = 1$. If we extend $\mathscr{L}$ to $\mathscr{L}^{+}$ by taking $\rho$ as a further initial function then either all the functions in $\mathscr{L}^{+}$ belong to $C^{\infty}$ or they are all continuous piecewise analytic functions. \newline
\indent Now set
\begin{equation} \label{6.1}
H(u,\vec{x}) = \rho(F(u,\vec{x})) \mbox{      }(\vec{x} = x_1,\ldots,x_k),
\end{equation}
and write $H_j(\vec{x})$ for $H(j,\vec{x})$. Then by 6.1 (3), (4), we have
\begin{equation} \label{6.2}
H_j(\vec{x}) = 0 \mbox{    }\textrm{for all }\vec{x}\textrm{ if }j \not\in A,
\end{equation}
\begin{equation} \label{6.3}
\exists\vec{x}\mbox{     }H_j(\vec{x}) = 1\textrm{ if }j \in A.
\end{equation}
But, by \eqref{3.2.6} and 6.1 (4) we see that, for $j \in A$,
\begin{equation} \label{6.4}
H_j(\vec{x}) = 0\textrm{ if }x_1^2,\ldots,x^2_k < \nu(j) - 1.
\end{equation}
Thus if an analogue machine is going to use $H_j$ to settle ?$j\in A$? and if $j \in A$, then the machine will have to calculate $P(j, y_1,\ldots,y_k)$ for some values $y_1,\ldots,y_k$ one at least of which - say $y_i$ - is greater than $\nu(j) - 1$. And by \eqref{3.2.7} the value of one of the $y$'s - $y_i$, say - must be accurate to within 1. Hence, for $j \in A$, the inputs $y_1,\ldots,y_k$ for the calculation of $H_j(y_1,\ldots,y_k)$ need to have a precision ratio of at least $\nu(j)$\footnote{Even if different PR's were used for $y_1,\ldots,y_k$ I believe the claim would stand: for the $m_1$ in \eqref{3.2.7} codes a computation sequence, so its size will certainly increase with $\nu(j)$.}. This is also true if $H_j$ is calculated by a digital computer. Thus the claim is proved for this example.
\subsection{}
Richardson, and following him, Da Costa and Doria make the problem look simpler by coding the $k$-plot $\vec{x}$ by a single real number $t$. Richardson defines decoding functions $(t)_1,\ldots,(t)_k$ (in $\mathscr{L}$) with the following property:\newline
\indent Given $\epsilon > 0$ and $x_1,\ldots,x_k$ one can find $t$ so that
\begin{equation}\label{6.3.1}
|x_i - (t)_i| < \epsilon\textrm{ for }1 \le i \le k.
\end{equation}
The functions he defines also satisfy
\begin{equation}\label{6.3.2}
(t)_i \le t.
\end{equation}
Now define a function $B_j$ by
\begin{equation} \label{6.3.3}
B_j(t) = H_j\left((t)_1,\ldots,(t)_k\right).
\end{equation}
Then
\begin{equation}\label{6.3.4}
B_j(t) = 0\textrm{ for all $t$, if }j\not\in A,
\end{equation}
while if $j \in A$ then for any $z < 1$
\begin{equation} \label{6.3.5}
\exists t\mbox{      }(B_j(t) > z).
\end{equation}
But, by \eqref{6.4} and \eqref{6.3.2} above we also have
\begin{equation} \label{6.3.6}
B_j(t) = 0\textrm{ if }t^2 < \nu(j) - 1.
\end{equation}
Any attempt to distinguish between \eqref{6.3.4} and \eqref{6.3.5} will yield further justifications for my claim. For example, Da Costa and Doria define
\begin{equation}\label{6.3.7}
K(j) = \int_0^{\infty}B_j(t)\gamma(t)dt
\end{equation}
where $\gamma(t)$ is a cut off factor inserted to ensure that the integral converges. (The exact nature of $B_j$ depends both on the distribution of the zeros of $P_A$ and on the particular decoding functions; in any case $B_j$ will be highly oscillatory, and, if $P_A$ has `rather few' zeros I think it likely that $\int_0^{\infty}B_j(t)dt$ will be of order ${\nu(j)}^{-1}$).\newline
\indent To specify and analogue machines which, for $j < J$ and $j \in A$ will output a non zero approximate value for $K(j)$ one will have to specify a value $B$ say, to replace $\infty$ as the upper limit of integration. But, by \eqref{6.3.5} above, one will then be able to compute a bound for $\beta(J)$ from $B$. And because of the cut off factor $\gamma$, \eqref{6.3.5} shows that $K(j)$ will be small of order ${\nu(j)}^{-1}$. Da Costa and Doria propose switching from one dynamical system to another, according to whether $K(j) = 0$ or $K(j) > 0$. An analogue machine which will correctly effect this switching will thus require, for the CVQ corresponding to $K(j)$ a precision ratio of order $\beta(J)$. Thus, in all, there are three different factors in the specification of the proposed machine which requires a knowledge of a bound for $\beta(J)$.

\section{Quantum Mechanical machines}
\subsection{}
Both my examples depend on specifying a self-adjoint operator $T$ on, say, Hilbert space (e.g. specifying the Hamiltonian for some quantum-mechanical system) and making observations on its spectrum to settle ?$j\in A$?.\newline
\indent The first example is due to Pour-El and Richards (CAP pp. 190-191). They show that a certain $T$ may be constructed as a computable limit of a sequence of computable operators $T_n$ with the following properties.\newline
\indent (1) Let $\lambda_j$ ($j \ge 0$) be a computable bounded sequence of real numbers. Then if $j \not\in A$ the spectrum of $T$ has $\lambda_j$ as an eigenvalue (corresponding to a line in spectranalytic terms), while if $j \in A$ the spectrum has a continuous band of width ${2.2}^{-\nu(j)}$ centered on $\lambda_j$. The factor $2^{-\nu(j)}$ ensures that the sequence $T_n$ has a computable modulus of convergence. To make observation easy one could take
\begin{equation*}
\lambda_j = 5-{4.2}^{-j},
\end{equation*}
and then there will be a gap between the bands (if present) around $\lambda_j$ and $\lambda_{j+1}$ to separate the lines or bands around $\lambda_j$ and $\lambda_{j+1}$ one only needs a precision of the order $2^{j}$; but to distinguish a \underline{line} at $\lambda_j$ and a \underline{band} around $\lambda_j$; one needs a precision of order $2^{\gamma(j)}$. Thus as in the previous examples, to settle ?$j \in A$? correctly for $j < J$ one needs to know a bound on $\beta(J)$ in order to ensure that the measyrements made will have the required precision. Another justification for my claim in this example is best illustrated by another example, which is a simplification of one given in Gandy (1991). Namely let the sequence $\{\lambda_n\}$ be defined by
\begin{equation*}
\lambda_n = 2^{-a(n)},
\end{equation*}
and let $S$ be a compact operator with these values of $\lambda_n$ as its eigenvalues. To decide ?$j\in A$? it is only necessary to observe, with say, a precision $2^{j+1}$, whether or not there is a line at $2^{-j}$. (Of course, on physical spectroscopy what one observes is transitions from one $\lambda$ to another, but this does not affect the argument.) So the question becomes: could one design a quantum mechanical device which would have, for some observable, an approximation $S'$ to $S$ whose eigenvalues for $j < J$ would be close to $S$? IT will be recalled that a design must allow one to compute approximate values for all relevant parameters and must specify allowed tolerances. I do not know, except in particular cases like atomic and molecular spectra, how one might construct a system which would approximate a given operator for a given observable. But it is obvious, for both $S$ and $T$, that one would need to know, at least approximately, the entries in the first $\beta(J)$ rows of their representing matrices (wrt some chosen orthonormal basis). But this justifies the claim\footnote{Both $S$ and $T$ are `effectively determined' operators. The interest of this concept lies not in examples like those given above but in the fact that the authors can (with considerable labour) give a general characterization, in terms of computability, for the spectra of such operators.}.
\subsection{} The wave functions for a quantum mechanical system may result from the superposition of infinitely many more easily defined wave functions and so correspond to the parallel working of infinitely many separate machines. This suggests a possible method for designing a quantum-mechanical device which would give correct answers to the questions ?$j \in A$? However the quantum computer described by Deutsch (1985) cannot do this, although it can use superposition greatly to reduce the run time for certain decidable problems.
\subsection{}
Refinements in experimental technique allow one to build analogue machines whose behaviour depends on a single quantum (e.g., a single photon). Experiments with such devices confirm the often counter-intuitive predictions of standard quantum theory. Could they provide a disproof of my claim? I do not know of any example for this.

\section{Discussion}
\subsection{}
When one shows that a given number-theoretic function is computable, or that a given number-theoretic problem is decidable, one does not place bounds on the run-time or the size of the memory - unless, of course, one is concerned with problems of complexity. That is, one is not concerned with precision ratios. So it may look as if I have placed unfair restrictions on analogue machines. But suppose one has proved that a certain programme will give correct answers to a problem ?$j \in X$?. Then, given $J$, one can \underline{compute} bounds on the time and space required to settle ?$j\in X$? correctly for all $j < J$. But this is exactly what I claim cannot be done for analogue machines intended to settle non-decidable problems.
\subsection{}
Cascades of events and chain reactions allow one (as in a photon multiplier) greatly to amplify the scale of an event. This is, in effect, a reduction of precision ratios. Could this be used to overcome the objections raised by my claim? The answer is `No', because only when one knows a bound for $\beta(J)$ can one determine how much amplification is needed.
\subsection{}
In CAP (and Pour-El \& Richards (1979)) other examples are given of differential equiations (in particular the wave equations) which will give a non-computable output for a computable input. The claim can be justified for these using the ideas of \S 5.
\subsection{}
Kreisel has discussed calculation by analogue machines in a number of place; see, in particular, his (1974), (1982), and (199 ). Some of his comments and analysis are illuminating, and have helped me in getting my ideas stragith. But one of his points is that there are more interesting, more sensible, and more relevant questions to ask than the (logical) question with which I am concerned.
\subsection{}
Penrose, in his (1989) and (1994), has argued that the human brain can be thought of as an analogue machine which can, in principle, settle undecidable problems. Firstly, he believes that mathematical results which can, at least in principle, be produced by human intelligence, cannot, even in principle, be produced by artificial intelligence - that is by some fixed programme $P$. Note that $P$ need not be itself directly responsible for the mathematical statements which the machine outputs. $P$ may be like an operating system, for example it may, by a process similar to natural selection, use mutations and tests of fitness to direct the (continual) evolution of subprogrammes for doing mathematics. But this possibility does not, straightforwardly, invalidate Penrose's argument justifying his belief. A concise version of Penrose's argument is given in Gandy (1994). Secondly Penrose believes that the sentences uttered or written by people are caused by physical and chemical events in their brains.\newline
\indent To allow for non-algorithmic actions in the brain, Penrose postulates a - not yet completely formulate - future theory which he calls CQG (for Correct Quantum Gravity). This will have consequences both for cosmology (concerning the direction of time's arrow) and for quantum theory (accounting for the collapse of real (not subjective) wave functions). He suggests ways in which such a theory may allow for the growth of microscopic structures (such as quasi-crystals, synapses and micro tubules in neurons) in ways which are not locally determined nor computable. It seems worthwhile to consider (rather naively) such patterns of grwoth from a mathematical point of view.

\section{Patterns of growth}
I consider a pattern of possible growth as being displayed on a tree. At each node $P$ there is a finite label which represents a particular structure $S_P$ at a particular stage of growth - for example, a particular quasi-crystal. If this structure $S_P$ is capable of growth then there will be a finite number of nodes $P_1,\ldots,P_k$ immediately below $P$; each of the structures $S_{P_1},\ldots,S_{P_k}$ arises from $S_P$ by a single step of growth (for example, by the addition of a single molecule). Two distinct structures $S_P$ and $S_Q$ may, in one step, grow into the same structure. Hence a node may have two different immediate predecessors; these trees are not the same as those standardly used in recursion theory. A node $P$ and the corresponding structure $S_P$ are \underline{fertile} if there is an infinite path through $P$. If $P$ is not fertile then, however $S_P$ may grow, after a finite number of steps it will become a structure which can grown no more.\newline
\indent Now we suppose that the label representing any structure $S$ is (coded by) a finite sequence $u$ of 0's and 1's. We may suppose that the significant features of $S$ can be computed from $u$. An infinite path gives an infinite sequence $u_1,u_2,\ldots,$ of binary sequences. We define the \underline{growth function} $\gamma$ along the path by $\gamma(u_n) = u_{n+1}$. If the sequence is computable then so is $\gamma$; in particular there is a Turing machine $M$ which, when presented with $u_n$ on its tape, will eventually replace it by $u_{n+1}$. Now the action of $M$ is certainly locally determined; it will, for example, in general, inspect each of the digits in $u_n$. We shall say that $\gamma$ (and the infinite sequence) are \underline{potentially locally determined}.
\subsection{}
Suppose we are given a tree of structures and a growth function $\gamma$ which satisfies the following conditions:
\begin{enumerate}[(i)]
\item If $u$ codes a fertile structure $S$, then $\gamma(u)$ codes a fertile structure into which $S$ can grown in a single step.
\item The function $\gamma$ is not potentially locally determined.
\end{enumerate}
Then, starting from any fertile structure $S$ and iterating $\gamma$ will produce a non-computable infinite sequence of structures.\newline
\indent If one could examine, say, the first $J$ structures in this sequence one could compute the first $J$ values of some non-computable function. The precision ratio of observation has to be sufficiently large to enable one to determine the codes $u$ for these $J$ structures; it might well be a computable function of $J$.
\subsection{}
Since quasi-crystals have been observed which contain a very large number of molecules, Penrose suggests that their growth is not a matter of chance, but is governed by some - as yet unformulated - laws of non-local actions. If, further, the theory involved actions which were not even potentially locally determined, then it would allow analogue machines to produce non-recursive functions. One would not expect the theory to be totally deterministic; indeed it is plausible that there are at least two distinct infinite paths through any fertile point of the tree, and hence continuum many such. Although each path yields a non-computable function, one cannot use it to settle a specified undecidable problem.\newline
\indent But for the growth of microstructures in the brain, which determine how neurons behave and how they affect each other, one would expect that certain particular paths would be selected on would be permitted.
\subsection{}
The definition of `potentially locally determined' can be made quite general by considering, in place of the Turing machine $M$, any mechanism which satisfies the principles of Gandy (1980) - in particular, of course, the principle of `local causation'. And then one has a converse to 9.2 - if the growth function along an infinite path is potentially locally determined, then the sequence of structures along it is computable.
\subsection{}
It is well-known that there are binary trees whose nodes form a recursive set, which have infinite paths but no computable infinite paths; using this fact one can for example describe a finite set of tiles which can tile the whole plane, but only in a non-computable way (see Hanf (1974)). Using the notion of trial and error predicates (see Putnam (1965)) we can see how the lattermost infinite path, $\lambda$ say, might be grown. A node is specified by a finite binary sequence $u$ which describes (with 0 for `Left' and 1 for `Right') the path from the vertex leading to it, and we consider $u$ also as the structure starting at $u$. The size of this is just the length of $u$. Now we define a computable sequence $u_n$ of nodes on the tree as follows.
\begin{enumerate}[(i)]
\item $u_0 = ()$ (the vertex of the tree).
\item If $u_n$ is not terminal (has nodes of the tree below it) then
\begin{equation*}
u_{n+1} = u_n0
\end{equation*}
\item Suppose $u_n$ is terminal and has the form $v0$ or $v011\ldots 1$ then
\begin{equation*}
u_{n+1} = v1
\end{equation*}
\end{enumerate}
\indent\indent Since no node on $\lambda$ is terminal, none of the $u_n$ can lie on the right of $\lambda$. Below any node $v$ which lies to the left of $\lambda$ (e.g.; 10 if $\lambda(1) = 1$ $\lambda(2) = 1$) there can only be finitely many nodes of the tree (since $v$ cannot be fertile). Hence for some $n$ we must have a $u_n$ lying to the right of $v$. Thuse for any $J$ there will be an $n_J$ such that $u_{n_J} = \lambda(1), \lambda(2),\ldots,\lambda(J-1)$.
\subsection{}
At first sight it might look as if this process of trial and error growth could be accomodated in some reasonable physical theory. But this is an illusion; for not only is $n_J$ not computable from $J$, but there can be no computable bound on the lenghts of the sequences $u_n$ with $n < n_J$ which have to be explored before $u_{n_J}$ is arrived at. And so the process considered is analogous to a trial and error process for deciding if $j \in A$ (as in \S 3) - one simple looks ahead to see if, for some $n$, $a(n) = j$.
\subsection{}
Penrose suggests that in a theory of quantum gravity the process of growth would be represented by a superposition of wave functions each corresponding to a particular pattern of growth, and that the effect of gravity would be to collapse the wave function, so that only constituents corresponding to patterns of growth capable of producing large structures would survive. To picture this process on the binary tree let the \underline{potential size}, $\pi(v)$ of a node $v$ be the maximum length of all nodes $u$ extending (or lying below) $v$. If $v$ is fertile we set $\pi(v) = \infty$. Then the proposed theory would ensure that any permitted vertex would grow to some node of great size, though (in the simple form in which I stated it) it would not guarantee growth along an infinite path. It would well be that for a given $J$ there would be a $k_J$ such that any node of size greater than $k_J$ would agree with $\lambda$ at the first $J$ places. But this fact will not allow us to compute values of $\lambda$ from observations on large structures which have developed, \underline{unless} we know some (necessarily non-computable) bounds for $k_J$. If a theory of growth of the kind considered is to stand up against our claim it looks as if some kind of non-computability must be built into the theory - for example into the way in which gravity determines the collapse of wave functions.
\newpage
\bibliographystyle{amsplain}

\printaddress
\makeatletter
\let\@setaddresses\relax
\makeatother

\includepdf[pages=-, addtotoc={1, section, 1, Gandy Original Handwritten Manuscript, GandyOriginal}]{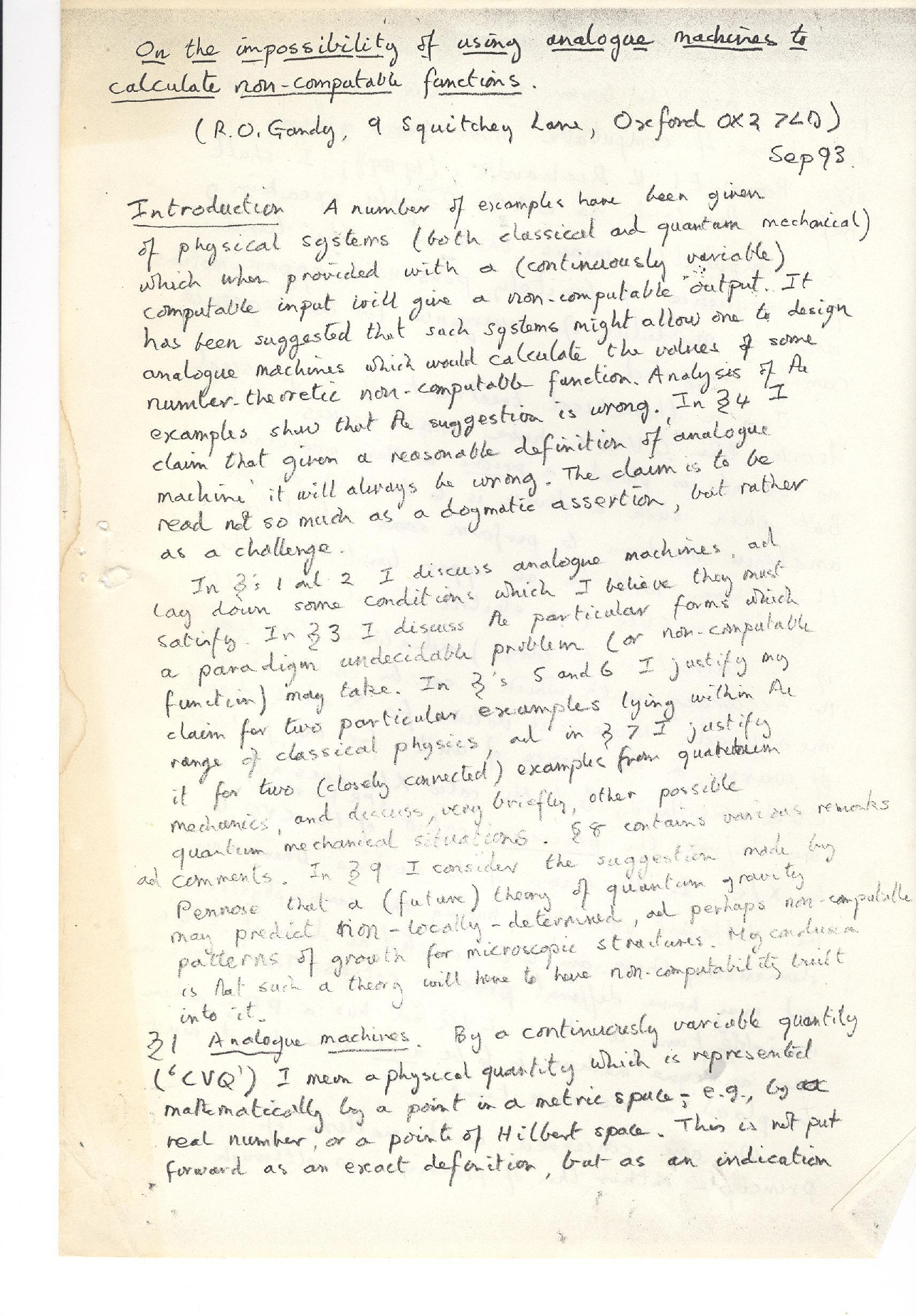}
\end{document}